\documentclass[12pt]{article}
\usepackage{bbm}
\usepackage{amsfonts}
\usepackage{pifont}
\usepackage{hyperref}
\usepackage{mathrsfs}
\usepackage{indentfirst}
\usepackage{amsmath}
\usepackage{amssymb}
\usepackage[amsmath, thmmarks]{ntheorem}
\usepackage{cite}
\usepackage{graphicx}
\usepackage{tikz}
\newtheorem{lemma}{\bf Lemma}[section]
\newtheorem{theorem}{\bf Theorem}[section]
\newtheorem{remark}{\bf Remark}[section]
\newtheorem{corollary}{\bf Corollary}[section]
\newtheorem{proposition}{\bf Proposition}[section]
\newtheorem{example}{\bf Example}[section]

\def\QEDopen{{\setlength{\fboxsep}{0pt}\setlength{\fboxrule}{0.2pt}\fbox{\rule[0pt]{0pt}{1.3ex}\rule[0pt]{1.3ex}{0pt}}}} 
\def\QED{\QEDopen}
\def\proof{{\bf Proof.} }
\def\endproof{\hspace*{\fill}~\QED\par\endtrivlist\unskip}

\begin{document}
\setcounter{page}{1}

\title{{\textbf{Conditions that the quotient of $L$-fuzzy up-sets forms a complete lattice}}\thanks {Supported by National Natural Science
Foundation of China (No. 61573240)}}
\author{Peng He\footnote{\emph{E-mail address}: 443966297@qq.com}, Xue-ping Wang\footnote{Corresponding author. xpwang1@hotmail.com; fax: +86-28-84761393}\\
\emph{College of Mathematics and Software
Science, Sichuan Normal University,}\\
\emph{Chengdu, Sichuan 610066, People's Republic of China}}

\newcommand{\pp}[2]{\frac{\partial #1}{\partial #2}}
\date{}
\maketitle
\begin{quote}
{\bf Abstract}

This paper deals with conditions under which the quotient of $L$-fuzzy up-sets forms a complete lattice by using terminologies of closure operators. It first gives a condition that a family of some subsets of a nonempty set can be represented by $L$-fuzzy up-sets, which is then used to formulate a
necessary and sufficient condition under which the quotient of $L$-fuzzy up-sets forms a complete lattice. This paper finally shows that the quotient of a kind of $L$-fuzzy up-sets is isomorphic to an interval generated by an $L$-fuzzy up-set.\\


\end{quote}
\section{Introduction}\label{intro}

$L$-fuzzy sets and structures have been widely studied from Goguen's first paper \cite{Goguen}. These structures appear
when the membership grades can be represented by elements of an ordered set, instead of just by numbers in the unit $[0,1]$. It is well-known that $L$-fuzzy mathematics attracts more and more interest in many
branches, for instance, algebraic theories including order-theoretic structures (see e.g., \cite{Seselja5}), automata and tree series
(see e.g., \cite{Borchar}) and theoretical computer science (see e.g., \cite{Chechik}). Among all the tropics on $L$-fuzzy mathematics, the representation of a poset by an $L$-fuzzy set is very interesting, which had been studied by Gorjanac-Rantiovi\'{c} and Tepav\v{c}evi\'{c} \cite{Gorjan} who formulated a necessary and sufficient condition, under which for a given family of subsets $\mathcal{F}$ of a set $X$ and a fixed complete lattice $L$ there is an $L$-fuzzy set $\mu$ such that the collection of cuts of $\mu$ coincides with $\mathcal{F}$. Jim\'{e}nez, Montes, \v{S}e\v{s}elja and Tepav\v{c}evi\'{c} \cite{Jorge} even showed a necessary and sufficient condition under which a collection of crisp up-sets (down-sets) of a poset $X$ consists of cuts of an $L$-fuzzy up-set (down-set). Moreover, He and Wang \cite{He} investigated the uniqueness of $L$-fuzzy sets whose collection of cuts is equal to a given family of subsets of a nonempty set.

In 2004, B. \v{S}e\v{s}elja and A. Tepav\v{c}evi\'{c} \cite{Seselja5} studied the $L$-fuzzy power set of a nonempty set, i.e., the collection of all $L$-fuzzy sets, and gave a necessary and sufficient condition under which two $L$-fuzzy sets from the collection have equal families of cuts. In this way, they obtained an equivalence relation on the $L$-fuzzy power set, and proved that the classes under the equivalence relation (the quotient of $L$-fuzzy sets) can be naturally ordered. They further provided a condition that the quotient of $L$-fuzzy sets forms a complete lattice.

Due to the importance of investigation of $L$-fuzzy up-sets and its equivalence, a natural problem is whether there
exists a similar result for $L$-fuzzy up-sets, that is to say, whether there exists a necessary and sufficient condition under which the quotient of $L$-fuzzy up-sets forms a complete lattice.
This paper will focus on the problem.

The paper is organized as follows. For the sake of convenience, some
notions and previous results are given in Section 2. In Section3, a condition that a family of some subsets of a nonempty set can be represented by $L$-fuzzy up-sets is obtained. The condition is further used to show a
necessary and sufficient condition under which the quotient of $L$-fuzzy up-sets forms a complete lattice, and it is verified that the quotient of a kind of $L$-fuzzy up-sets is isomorphic to an interval generated by an $L$-fuzzy up-set in Section 4. Conclusions are drawn
in Section 5.

\section{Preliminaries}
In this section, we first list some necessary notions and relevant properties from the classical order theory in the sequel. For more comprehensive presentation, see e.g., book \cite{Gratze}.

A poset is a structure $(P, \leq)$ (or $P$ for short) where $P$ is a nonempty set and $\leq$ an ordering (reflexive, antisymmetric and transitive)
relation on $P$. A sub-poset of a poset $(P, \leq)$ is a poset $(Q, \leq)$ where $Q$ is a nonempty subset of $P$ and $\leq$ on $Q$ is
restricted from $P$. A lattice is a poset $L$ in which for each pair of elements $x, y$ there is a greatest lower bound (glb, infimum, meet) and a least upper bound (lub, supremum, join), denoted, respectively, by $x\wedge y$ and $x\vee y$. These are binary operations on $L$. A poset $L$ is said to be a complete lattice if infima and suprema exist for each subset of $L$. A complete lattice $L$ possesses the top element
$1_L$ and the bottom element $0_L$. An up-set (semi-filter) on a poset $P$ is any sub-poset $U$, satisfying: for $x\in U, y\in P, x\leq y$
implies $y\in U$. In particular, for an $x\in P$, a sub-poset $\uparrow x:=\{y\in P\mid x\leq y\}$ is a principle filter generated by $x$.

We say that an element $q$ in a lattice $L$ is meet-irreducible if, $q\neq 1_L$ and for all $x, y\in L$, $q=x\wedge y$ implies $q=x$ or $q=y$.

A closure operator on a poset $P$ is a function $C: P\rightarrow P$ such that , for all $p, q\in P$, it fulfills
the following three requirements:\\
(a) $p\leq C(p)$,\\
(b) $p\leq q\Rightarrow C(p)\leq C(q)$,\\
(c) $C(C(p))=C(p)$.\\
If $p=C(p)$, then $p$ is a closed element under the corresponding closure operator.


\begin{lemma}[\cite{Gratze}]\label{le 2}
Let $L$ be a complete lattice. Then, for any closure operator $C$ on $L$, the set of all closed elements of $L$ is a
complete lattice under the order inherited from $L$.
\end{lemma}
\begin{lemma}[\cite{Gratze, Jorge}]\label{le 3}
Let $C$ be a closure operator on a poset $L$. If we consider the relation $\sim'$ on $L$ defined by $x\sim' y$ if and only if $C(x)=C(y)$, then:\\
\emph{(1)} $\sim'$ is an equivalence relation on $L$.\\
\emph{(2)} Each equivalence class possesses the top element which is closed under $C$.\\
\emph{(3)} The set $L/{\sim'}$ can be ordered: $[x]\leq [y]$ if and only if $C(x)\leq C(y)$ in $L$, where $[x]=\{a\in L\mid x\sim' a\}$ for all $x\in L$.

Furthermore, if $L$ is a complete lattice, then: \\
\emph{(4)} The poset $(L/{\sim'}, \leq)$ is a lattice isomorphic with the poset of all closed elements of $L$ under $C$.
\end{lemma}

In what follows, we denote $L/{\sim'}$ by $L/ C$ if $C$ is a closure operator on a poset or a complete
lattice $L$. Let $C(L)=\{C(a)\mid a\in L\}$. Then from Lemma \ref{le 2}, we know that the poset $C(L)$ is a complete lattice if $L$ is a complete lattice. Further, we have the following two propositions.
\begin{proposition}\label{pro 1}
Let $C_0$ be a closure operator on a complete lattice $L$ and let $C_1$ be a closure operator on the complete
lattice $C_0(L)$. Then the function $C: L\rightarrow L$ associated
to any element $p\in L$ the value $C_1(C_0(p))$ is a closure operator on $L$, and we denote $C=C_1 \circ C_0$.
\end{proposition}
\proof We need to prove that $C$ satisfies the three conditions of the definition of a closure operator.

From the hypotheses of Proposition \ref{pro 1}, $p\leq C_0(p)$ and $C_0(p)\leq C_1(C_0(p))$ for all $p\in L$. Then $p\leq C(p)$. If $p\leq q$ for all $p, q\in L$, then $C_0(p)\leq C_0(q)$ since $C_0$ is a closure operator, and $C(p)=C_1(C_0(p))\leq C_1(C_0(q))=C(q)$ since $C_1$ is also a closure operator.

Note that $C_1(C_0(p))\in C_0(L)$ for all $p\in L$. Then $$C_0(C_1(C_0(p)))=C_1(C_0(p)).$$
Thus, $C_1\circ C_0(C_1(C_0(p)))=C_1(C_1(C_0(p))$. Since $C_1$ is a closure operator on $C_0(L)$ and $C_1(C_0(p))\in C_1(C_0(L))$, we have
$$C_1(C_1(C_0(p))=C_1(C_0(p)).$$ The above two equalities imply that $C_1\circ C_0(C_1(C_0(p)))=C_1(C_0(p))$, i.e., $C(C(p))=C(p)$.

We conclude that $C=C_1\circ C_0$ is a closure operator on $L$.
\endproof

\begin{proposition}\label{pro 2}
Under the assumptions of Proposition \ref{pro 1}, $L/(C_1\circ C_0) \cong (L/C_0)/C_1$.
\end{proposition}
\proof
By Lemma \ref{le 3}, $L/C_0 \cong C_0(L)$ and $C_0(L)/ C_1\cong C_1\circ C_0 (L)$.
Thus $$(L/C_0)/C_1 \cong C_0(L)/C_1 \cong C_1\circ C_0 (L).$$
Moreover, by Proposition \ref{pro 1}, we know that $C_1\circ C_0$ is a closure operator on $L$.
Then from Lemma \ref{le 3}, $$L/(C_1\circ C_0)\cong C_1\circ C_0(L).$$

Therefore, $L/(C_1\circ C_0) \cong (L/C_0)/C_1$.
\endproof

 In the following, we shall recall some notions of $L$-fuzzy sets. More details about the relevant properties can be found e.g.,
in \cite{Goguen,Gorjan, Jorge, Seselja5}.

An $L$-fuzzy set, Lattice-valued or $L$-valued set here is a mapping $\mu : X\rightarrow L$ from a nonempty set $X$ into a complete
lattice $(L, \leq)$. Let $\mathcal{F}_{L}(X)$ be the collection of all $L$-fuzzy sets on $X$.

If $\mu : X\rightarrow L$ is an $L$-fuzzy set on $X$ then, for $p\in L$, the set \begin{equation}\label{equ 00}\mu_p =\{x\in X\mid \mu(x)\geq p\}\end{equation} is called the $p$-cut,
a cut set or simply a cut of $\mu$. Let $\mu_L=\{\mu_p\mid p\in L\}$ for all $\mu\in \mathcal{F}_{L}(X)$.

An $L$-fuzzy up-set is a mapping $\mu : X\rightarrow L$ from a poset $(X, \leq)$ into a complete
lattice $(L, \leq)$ such that for all $x, y\in X$ $$x\leq y\Rightarrow \mu(x)\leq \mu(y).$$
Let $\mathcal{F}^u_{L}(X)$ be the collection of all $L$-fuzzy up-sets on $X$. Obviously, every $L$-fuzzy up-set on $X$ is also an $L$-fuzzy set.

Let $\mathcal{P}(X)$ be the power set of a set $X$. Then the following two statements characterize the collection of cuts of an $L$-fuzzy up-set.
\begin{theorem}[\cite{Jorge}]\label{the 1}
Let $(X, \leq)$ be a poset, let $\mathcal{F}\subseteq\mathcal{P}(X)$ be a family of some up-sets of a poset $X$, and let $(L, \leq)$ be a complete lattice.
Then, there is an $L$-fuzzy up-set $\mu: X\rightarrow L$ such that its family of cuts is equal to $\mathcal{F}$ if and only if the following two conditions hold:\\
(1) $\mathcal{F}$ is closed under intersections and contains $X$;\\
(2) there is a closure operator $C$ on $L$, such that the poset $(L/C, \leq)$ is order isomorphic to $(\mathcal{F}, \supseteq)$.
\end{theorem}
\begin{theorem}[\cite{Jorge}]\label{the 2}
Let $(X, \leq)$ be a poset, let $L$ be a complete lattice and let $\mu$ be an L-fuzzy set on $X$. Then the following two statements are equivalent: \\
(1) $\mu$ is an L-fuzzy up-set on $X$.\\
(2) The p-cut $\mu_p$ of $\mu$ is an up-set on $X$, for any $p\in L$.
\end{theorem}

\section{Representation of families of subsets by $L$-fuzzy up-sets}

In this section, we shall investigate a condition that a family of some subsets of a nonempty set can be represented by $L$-fuzzy up-sets.

Let $A$ and $B$ be two sets.
Then we denote that $A\setminus B=\{x\in A\mid x\notin B\}$, for convenience, if $B=\{b\}$ then we
write $A\setminus B$ as $A\setminus b$. We denote by $M(L)$ the set of all meet-irreducible elements of a lattice $L$.

\begin{theorem}\label{the0 1}
Let $(X, \leq)$ be a poset, $(L, \leq)$ a complete lattice and $\mu \in \mathcal{F}^{u}_{L}(X)$.
If $\mathcal{T}\subseteq \mu_L$ is closed under intersections and contains $X$, then there
exists $\nu\in \mathcal{F}^{u}_{L}(X)$ such that $\nu_L=\mathcal{T}$.
\end{theorem}
\proof
Suppose that $F=(\mu_L, \supseteq)$ and $T=(\mathcal{T}, \supseteq)$. Then both $F$ and $T$ are complete lattices.
Let $C_1: F\rightarrow F$ be defined by
\begin{equation}\label{equa 1}
C_1(p)=\left\{
\begin{array}{rcl}
\overline{\bigcup}\{q\in \mathcal{T}\setminus 1_T \mid q\subsetneq p \}       &      & {p \in \mathcal{T}\setminus M(T),}\\
$ $\\
\overline{\bigcup}\{q\in \mathcal{T}\setminus 1_T \mid q\subseteq p \}     &      & {\mbox{ otherwise } }
\end{array}\right.
\end{equation}
in which, $\overline{\bigcup} S=\bigwedge_T S$ when $S\neq \emptyset$, and $\overline{\bigcup} \emptyset =1_F$ where $\bigwedge_T$ is the meet of $T$.

By (\ref{equa 1}), we have that \begin{equation}\label{equa 2}C_1(F)\subseteq \mathcal{T} \bigcup \{1_F\}\end{equation} since $T$ is a complete lattice.
Moreover, we claim that \begin{equation}\label{equa 6}1_T\notin C_1(F)\end{equation} when $1_T\neq 1_F$.
Otherwise, there exists $p\in \mu_L$ such that $C_1(p)=1_T$, which means that $\overline{\bigcup}\{q\in \mathcal{T}\setminus 1_T \mid q\subsetneq p\}=1_T$ or
$\overline{\bigcup}\{q\in \mathcal{T}\setminus 1_T \mid q\subseteq p\}=1_T$ by (\ref{equa 1}).
Suppose $\overline{\bigcup}\{q\in \mathcal{T}\setminus 1_T \mid q\subseteq p\}=1_T$. Clearly, $\{q\in \mathcal{T}\setminus 1_T \mid q\subseteq p\}\neq \emptyset$
since $\overline{\bigcup} \emptyset=1_F\neq 1_T$. Then there exists an element $q\in \mathcal{T}\setminus 1_T$ such that $q\subseteq p$, which together with formula (\ref{equa 1}) implies
that $C_1(p)\supseteq q\supsetneq 1_T$, a contradiction. Similarly, the equality $\overline{\bigcup}\{q\in \mathcal{T}\setminus 1_T \mid q\subsetneq p\}=1_T$
will deduce a contradiction.

The following proof is made in three steps:

First, we prove $C_1$ is a closure operator on $F$.

It suffices to prove that $C_1$ satisfies the three conditions (a), (b) and (c). Obviously, by formula (\ref{equa 1}), $p\supseteq C_1(p)$ for all $p\in \mu_L$, i.e., the condition (a) holds.

Assume that $p, r\in \mu_L$ and $p\supseteq r$. It is clear that $p=r$ implies $C_1(p)=C_1(r)$. Next, we suppose $p\supsetneq r$.
Then from (\ref{equa 1}), $$C_1(r)\subseteq\overline{\bigcup}\{q\in \mathcal{T}\setminus 1_T \mid q\subseteq r\} \mbox{ and } C_1(p)\supseteq\overline{\bigcup}\{q\in \mathcal{T}\setminus 1_T \mid q\subsetneq p\}.$$
As $p\supsetneq r$, we have $$C_1(p)\supseteq\overline{\bigcup}\{q\in \mathcal{T}\setminus 1_T \mid q\subsetneq p\}\supseteq\overline{\bigcup}\{q\in \mathcal{T}\setminus 1_T \mid q\subseteq r\}\supseteq C_1(r).$$
Therefore, $C_1(p)\supseteq C_1(r)$, i.e., the condition (b) holds.

Next we prove the condition (c) is true.

Note that $C_1(p)\supseteq C_1(C_1(p))$ for all $p\in \mu_L$
since the condition (a) holds. Thus it suffices to prove that $$C_1(C_1(p)) \supseteq C_1(p) \mbox{ for all }p\in \mu_L.$$

By formula (\ref{equa 2}), there are three cases as follows.

Case (1). If $C_1(p)=1_F$, then from (\ref{equa 1}) \begin{equation}\label{equa 3}C_1(C_1(p))=C_1(1_F)=\overline{\bigcup} \emptyset =1_F=C_1(p).\end{equation}

Case (2). If $C_1(p)\in M(T)$ then $$C_1(C_1(p))=\overline{\bigcup}\{q\in \mathcal{T}\setminus 1_T \mid q\subseteq C_1(p)\}\supseteq C_1(p)$$
since $C_1(p)\in \{q\in \mathcal{T}\setminus 1_T \mid q\subseteq C_1(p)\}$.

Case (3). If $C_1(p)\in (\mathcal{T}\setminus M(T))\setminus 1_F$, then we claim that $C_1(p)\neq 1_T$. Otherwise, $C_1(p)=1_T\in C_1(F)$,
which means that $1_T=1_F$ by (\ref{equa 6}). Thus $C_1(p)=1_F\in (\mathcal{T}\setminus M(T))\setminus 1_F$, a contradiction. Hence $C_1(p)\neq 1_T$. Then
there exist two element $m, n\in \mathcal{T}\setminus 1_T$ such that $m\subsetneq C_1(p)$, $n\subsetneq C_1(p)$ and $m\wedge_T n=C_1(p)$.
Consequently, $$C_1(C_1(p))=\overline{\bigcup}\{q\in \mathcal{T}\setminus 1_T \mid q\subsetneq C_1(p)\}\supseteq m\overline{\cup}n= m\wedge_T n=C_1(p).$$

Cases (1), (2) and (3) imply that $C_1(C_1(p))\supseteq C_1(p)$ for all $p\in \mu_L$.

According to the above proof, we conclude that $C_1$ is a closure operator on $F$.

Secondly, we shall prove \begin{equation}\label{eq111}F/C_1 \cong T.\end{equation}

Let $\mathcal{S}=(\mathcal{T}\setminus 1_T) \bigcup \{1_F\}$. Then we easily know that $(\mathcal{S}, \supseteq)\cong T$. On the other hand, by Lemma \ref{le 3}, we have $F/C_1\cong (C_1(F), \supseteq)$ since $C_1$ is a closure operator on $F$. Thus, in order to prove $F/C_1 \cong T$, we just need to prove $\mathcal{S}=C_1(F)$.

By (\ref{equa 2}) and (\ref{equa 6}), we have that $C_1(F) \subseteq \mathcal{S}$. Now we shall prove $C_1(F) \supseteq \mathcal{S}$.

From formula (\ref{equa 3}), $C_1(1_F)=1_F$. By the construction of $\mathcal{S}$, we have that $x\in M(T)$ or $x\in (\mathcal{T}\setminus M(T))\setminus 1_T$
for all $x\in \mathcal{S}$ with $x\neq 1_F$. If $x\in M(T)$ then $x\in \{q\in \mathcal{T}\setminus 1_T \mid q\subseteq x\}$, which yields
that $C_1(x)=\overline{\bigcup}\{q\in \mathcal{T}\setminus 1_T \mid q\subseteq x\}\supseteq x$. Note that $x\supseteq C_1(x)$ since $C_1$ is a closure operator on $F$. Thus $C_1(x)=x$.
If $x\in (\mathcal{T}\setminus M(T))\setminus 1_T$ then there exist two element $y, z\in \mathcal{T}\setminus 1_T$ with $y\subsetneq x$ and $z\subsetneq x$ such that $y\wedge_T z=x$.
Thus $$C_1(x)=\overline{\bigcup}\{q\in \mathcal{T}\setminus 1_T \mid q\subsetneq x\}\supseteq y\overline{\cup} z=y\wedge_T z= x.$$  Then $C_1(x)=x$ since $x\supseteq C_1(x)$. Therefore, we always have that $C_1(x)=x$ for all $x\in \mathcal{S}$ with $x\neq 1_F$, i.e., $C_1(F) \supseteq \mathcal{S}$.

Finally, we shall prove that there exists $\nu\in \mathcal{F}^{u}_{L}(X)$ such that $\nu_L=\mathcal{T}$.

In fact, by Theorem \ref{the 1} and the hypotheses of Theorem \ref{the0 1}, there exists a closure operator $C_0$ on $L$ such that $L/ C_0\cong F$ since $\mu_L=\mu_L\subseteq\mathcal{P}(X)$. Then by Lemma \ref{le 3},
\begin{equation}\label{121}L/ C_0\cong C_0(L)\cong F.\end{equation} From formulas (\ref{121}) and (\ref{eq111}), we easily prove that there exists a closure
operator $C_*$ on $C_0(L)$ such that $C_0(L)/C_*\cong T$. Let $C=C_0\circ C_*$. Then by Propositions \ref{pro 1} and \ref{pro 2}, $C$ is a
closure operator on $L$ and $L/C\cong (L/C_0)/C_*\cong T$. Again by Theorem \ref{the 1}, there exists $\nu\in \mathcal{F}^{u}_{L}(X)$ such that $\nu_L=\mathcal{T}$.

This completes the proof.
\endproof

The following example illustrates an application of Theorem \ref{the0 1}.
\begin{example}\label{exa 00}
\emph{Let us consider the poset $(X, \leq)$ and the complete lattice $(L, \leq)$ represented in Fig. 1. Let
$\mathcal{R}=\{\emptyset, \{a\}, \{b\}, \{b, c\}, \{a, b, c, d, e\}\}$ be a family of some up-sets of $X$. Does there exist $\mu \in\mathcal{F}^u_{L}(X)$ such that $\mu_L=\mathcal{R}$ ?}
\par\noindent\vskip50pt
 \begin{minipage}{11pc}
\setlength{\unitlength}{0.75pt}\begin{picture}(600,100)
\put(330,40){\circle{4}}\put(335,30){\makebox(0,0)[l]{\footnotesize $0_L$}}
\put(295,75){\circle{4}}\put(290,65){\makebox(0,0)[l]{\footnotesize $q$}}
\put(330,75){\circle{4}}\put(335,70){\makebox(0,0)[l]{\footnotesize $r$}}
\put(365,75){\circle{4}}\put(370,70){\makebox(0,0)[l]{\footnotesize $p$}}
\put(400,75){\circle{4}}
\put(400,110){\circle{4}}
\put(347.5,127.5){\circle{4}}
\put(312.5,92.5){\circle{4}}
\put(295,110){\circle{4}}\put(280,110){\makebox(0,0)[l]{\footnotesize $s$}}
\put(365,110){\circle{4}}\put(370,110){\makebox(0,0)[l]{\footnotesize $t$}}
\put(330,145){\circle{4}}\put(335,150){\makebox(0,0)[l]{\footnotesize $1_L$}}
\put(295,77){\line(0,1){31}}
\put(328.5,41.5){\line(-1,1){32}}
\put(331.5,41.5){\line(1,1){32}}
\put(331.5,41.5){\line(2,1){66.5}}
\put(400,77){\line(0,1){31}}
\put(398.5,111.5){\line(-3,1){49.5}}
\put(398.5,76.5){\line(-1,1){32}}
\put(330,42){\line(0,1){31}}
\put(331.5,76.5){\line(1,1){32}}
\put(328.5,76.5){\line(-1,1){14.5}}
\put(311,94){\line(-1,1){14.5}}
\put(314,94){\line(1,1){32}}
\put(365,77){\line(0,1){31}}
\put(296.5,111.5){\line(1,1){32}}
\put(363.5,111.5){\line(-1,1){14.5}}
\put(346,129){\line(-1,1){14.5}}
\put(330,5){$L$}
\put(165,40){\circle{4}}\put(160,30){\makebox(0,0)[l]{\footnotesize $e$}}
\put(235,40){\circle{4}}\put(240,30){\makebox(0,0)[l]{\footnotesize $c$}}
\put(200,75){\circle{4}}\put(205,80){\makebox(0,0)[l]{\footnotesize $b$}}
\put(200,40){\circle{4}}\put(205,30){\makebox(0,0)[l]{\footnotesize $d$}}
\put(165,75){\circle{4}}\put(155,80){\makebox(0,0)[l]{\footnotesize $a$}}
\put(166.5,41.5){\line(1,1){32}}
\put(165,42){\line(0,1){31}}
\put(198.5,41.5){\line(-1,1){32}}
\put(233.5,41.5){\line(-1,1){32}}
\put(195,5){$X$}
\put(155,-15){\emph{Fig.1 Hasse diagrams of $X$ and $L$}}
  \end{picture}
  \end{minipage}

$$ $$

\emph{Let $\mathcal{S}=\{\emptyset, \{a\}, \{b\}, \{a, b\}, \{b, c\}, \{a, d\}, \{a, b, c, d, e\}\}$. Then one can check that $\mathcal{S}$ is also a family of some up-sets of $X$ and there exists a closure operator $C$ on $L$ such that $$L/C\cong (\{0_L, p, q, r, s, t, 1_L\}, \leq)\cong (\mathcal{S}, \supseteq).$$ Thus, by Theorem \ref{the 1}, there exists $\nu\in\mathcal{F}^u_{L}(X)$ such that $\nu_L=\mathcal{S}$.
Since $\mathcal{R}\subseteq \mathcal{S}$, there exists $\mu \in\mathcal{F}^u_{L}(X)$ such that $\mu_L=\mathcal{R}$ by Theorem \ref{the0 1}.}
\end{example}

The following remark is clear.
\begin{remark}\label{re}\emph{Each nonempty set $X$ can be regarded as a special poset for which if $x, y\in X$ with $x\neq y$ then $x$ and $y$ are not comparable.}\end{remark}

Then from Remark \ref{re} and Theorem \ref{the0 1}, we have the following corollary.
\begin{corollary}\label{cor0 1}
Let $X$ be a nonempty set, $(L, \leq)$ a complete lattice and $\mu \in \mathcal{F}_{L}(X)$.
If $\mathcal{T}\subseteq \mu_L$ is closed under intersections and contains $X$, then there
exists $\nu\in \mathcal{F}_{L}(X)$ such that $\nu_L=\mathcal{T}$.
\end{corollary}
\section{Conditions that $(\mathcal{F}^{u}_{L}(X)/\sim, \leq)$ is a complete lattice}
In this section, we shall show necessary and sufficient conditions that $(\mathcal{F}^{u}_{L}(X)/\sim, \leq)$ is a complete lattice, and prove that $(\mathcal{F}^{u}_{L}(X/C)/\sim, \leq)$ is isomorphic to an interval generated by an $L$-fuzzy up-set, in which $C$ is a closure operator on the poset $X$.

Consider the relation $\sim$ on $\mathcal{F}_{L}(X)$ defined by $\mu \sim \nu$ if and only if $\mu_L=\nu_L$. It is easy to see that $\sim$ is an equivalence relation. Further, let us define the relation $\leq$ on $\mathcal{F}_{L}(X)/\sim$ in the following way:
$$[\mu]_{\sim}\leq [\nu]_{\sim} \mbox{ if and only if } \mu_L \subseteq \nu_L.$$

Given an $L$-fuzzy set $\mu\in \mathcal{F}_{L}(X)$, we define the relation $\approx_\mu$ on $L$ by $p\approx_\mu q$ if and only if $\mu_p=\mu_q$. Then $\approx_\mu$ is obviously an equivalence relation on $L$. In 2004, \v{S}e\v{s}elja and Tepav\v{c}evi\'{c} \cite{Seselja5} proved the following statement.
\begin{proposition}\label{prop 1} Let $\mu\in \mathcal{F}_{L}(X)$. Then
$(L/\approx_\mu, \leq)\cong (\mu_L, \supseteq)$.
\end{proposition}

They further gave a condtion that $(\mathcal{F}_{L}(X)/\sim, \leq)$ is a complete lattice as follows.
\begin{theorem}[\cite{Seselja5}]\label{the 4}
If there exists $\mu \in \mathcal{F}_{L}(X)$ such that $L/\approx_\mu$ is isomorphic to the power set of a set $X$
 then $(\mathcal{F}_{L}(X)/\sim, \leq)$ is a complete lattice.
\end{theorem}

In what follows, we shall give conditions that $(\mathcal{F}^{u}_{L}(X)/\sim, \leq)$ is a complete lattice.

Given any set $A$, we denote its cardinality by $|A|$. We first have the following statement.
\begin{proposition}\label{pro 3}
Let $(X, \leq)$ be a poset and $|L|=1$. Then $(\mathcal{F}^{u}_{L}(X)/\sim, \leq)$ is a complete lattice.
\end{proposition}
\proof
Clearly, $\mu_L=\{X\}$ for all $\mu\in \mathcal{F}^u_{L}(X)$. Therefore, $(\mathcal{F}^u_{L}(X)/\sim, \leq)$ is surely a complete lattice.
\endproof

From Remark \ref{re} and Proposition \ref{pro 3}, the following corollary is immediately.
\begin{corollary}\label{cor 3}
Let $X$ be a nonempty set and $|L|=1$. Then $(\mathcal{F}_{L}(X)/\sim, \leq)$ is a complete lattice.
\end{corollary}

Let $(X, \leq)$ be a poset and let $\mathcal{F}_X$ be the set of all up-sets of $X$. We denote $(\mathcal{F}_X, \supseteq)$ by $F_X$. Then we have the following theorem.
\begin{theorem}\label{the0 2}
Let $(X, \leq)$ be a poset and $L$ a complete lattice with $|L|>1$. Then $(\mathcal{F}^{u}_{L}(X)/\sim, \leq)$ is a complete lattice if and only if there exists a closure operator $C$ on $L$ such that $L/C\cong F_X$.
\end{theorem}
\proof
Necessity. Suppose that $(\mathcal{F}^{u}_{L}(X)/\sim, \leq)$ is a complete lattice. We first note that $0_L\neq 1_L$ since $|L|>1$, and there exists a closure operator $C$ on $L$ such
that $C(L)=\{0_L, 1_L\}$. Then, by Lemma \ref{le 3}, $ L/C\cong (\mathcal{T},\supseteq)$ for all $\mathcal{T}=\{S, X\}$ where $S\subsetneq X$ and $S\in \mathcal{F}_X$.
Let $$\mathcal{S'}=\{\mathcal{T}\mid \mathcal{T}=\{S, X\}, S\subsetneq X \mbox{ and }S\in \mathcal{F}_X\}.$$
Then, by Theorem \ref{the 1},
there exists $\nu^{\mathcal{T}}\in \mathcal{F}^{u}_{L}(X)$ such that $\nu^{\mathcal{T}}_L=\mathcal{T}$ for all $\mathcal{T}\in \mathcal{S'}$.
Since $(\mathcal{F}^u_{L}(X)/\sim, \leq)$ is a complete lattice, we can let $\bigvee_{\mathcal{T}\in \mathcal{S'}}[\nu^{\mathcal{T}}]_{\sim}=[\mu]_{\sim}$. Then $[\mu]_{\sim}\in \mathcal{F}^u_{L}(X)/\sim$. Thus $\mu_L\supseteq \bigcup_{\mathcal{T}\in \mathcal{S'}}\nu^{\mathcal{T}}_L=\mathcal{F}_X$. Note that $\mu_L\subseteq \mathcal{F}_X$ by Theorem \ref{the 2}. Hence $\mu_L=\mathcal{F}_X$. Again, by Theorem \ref{the 1}, there exists
a closure operator $C$ on $L$ such that $L/C\cong F_X$.

Sufficiency. If there exists a closure operator $C$ on $L$ such that $L/C\cong F_X$, then, by Theorem \ref{the 1}, there exists $\mu \in \mathcal{F}^{u}_{L}(X)$ such that \begin{equation}\label{0011}\mu_L=\mathcal{F}_X.\end{equation}
Thus, by Theorem \ref{the 2}, we know that $[\mu]_\sim$ is the greatest element of $(\mathcal{F}^{u}_{L}(X)/\sim, \leq)$ since $\mathcal{F}_X$ is the set of all up-sets of $X$.
As it is well known, if a poset is closed under arbitrary infima and it has the greatest element, then the poset is a complete lattice.
Thus, in order to prove that $(\mathcal{F}^{u}_{L}(X)/\sim, \leq)$ is a complete lattice, it suffices to prove that $(\mathcal{F}^{u}_{L}(X)/\sim, \leq)$ is closed under arbitrary infima.

Let $\emptyset\neq B\subseteq \mathcal{F}^{u}_{L}(X)$. Then by Theorem \ref{the 1}, $X\in \nu_L$ and $\nu_L$ is closed under intersections since $\nu_L=\nu_L\subseteq\mathcal{P}(X)$ for all $\nu \in B$.
Suppose that $\mathcal{S}=\bigcap\{\nu_L\mid \nu\in B\}$. Then $\mathcal{S}$ is closed under intersections and contains $X$. Note that $\mathcal{S}\subseteq \mathcal{F}_X$.
Thus, using Theorem \ref{the0 1}, there exists $f\in \mathcal{F}^{u}_{L}(X)$ such that $f_L =\mathcal{S}$ since $\mathcal{F}_X=\mu_L$ by (\ref{0011}). Each collection of cuts
corresponding to a $\sim$-class, therefore, $(\mathcal{F}^{u}_{L}(X)/\sim, \leq)$ is closed under arbitrary infima.
\endproof

 Then by Remark \ref{re} and Theorem \ref{the0 2} we have:

\begin{corollary}\label{coro 1}
Let $X$ be a nonempty set and $L$ a complete lattice with $|L|>1$. Then $(\mathcal{F}_{L}(X)/\sim, \leq)$ is a complete lattice if and only if
there exists a closure operator $C$ on $L$ such that $L/C\cong(\mathcal{P}(X), \supseteq)$.
\end{corollary}

\begin{remark}\label{rem 1}
\emph{(i) Corollary \ref{coro 1} implies that the inverse of Theorem \ref{the 4} holds when $|L|>1$. In fact, by Corollary \ref{coro 1}, there exists a closure operator $C$ on $L$ such that $L/C\cong(\mathcal{P}(X), \supseteq)$. Since $\mathcal{P}(X)$ can be regarded as a family of all up-sets of $X$ by Remark \ref{re}, from Theorem \ref{the 1}, there exists $\mu \in \mathcal{F}^{u}_{L}(X)$ such that $\mu_L=\mathcal{P}(X)$. Furthermore, by Proposition \ref{prop 1} $(L/\approx_\mu, \leq)\cong (\mu_L, \supseteq)$. Then $(L/\approx_\mu, \leq)\cong (\mathcal{P}(X), \supseteq)$}.

\emph{(ii) If $|L|=1$ then the inverse of Theorem \ref{the 4} is not true. Indeed, if $|L|=1$, then by Corollary \ref{cor 3} we know that $(\mathcal{F}_{L}(X)/\sim, \leq)$ is a complete lattice. However, $(L/\approx_{\mu},\leq)\ncong (\mathcal{P}(X), \supseteq)$ for all $\mu\in\mathcal{F}_{L}(X)$ since $|\mathcal{P}(X)|>1$.}
\end{remark}

In what follows, we shall prove that $(\mathcal{F}^{u}_{L}(X/C)/\sim, \leq)$ is isomorphic to an interval generated by an $L$-fuzzy up-set, where $C$ is a closure operator on the poset $X$.

Note that it is easy to verify that $F_X$ is a complete distributive lattice for which the meet (join) is $\cup$ ($\cap$) (see \cite{Jorge}). Then we have the following theorem:

\begin{theorem}\label{the0 3}
Let $(X, \leq)$ be a poset and $C$ a closure operator on $X$. Then the lattice $F_{X/C}$ can
be embedded into $F_X$, such that all infima, suprema, and the top and bottom elements are preserved under the embedding.
\end{theorem}
\proof
Let \begin{equation}\label{equa 4}S_T=\bigcup_{[x]\in T}\{y\in X\mid y\in [x]\}\end{equation} for all $T\in \mathcal{F}_{X/C}$ and
\begin{equation}\label{equa 5}\mathcal{S}=\{S_T\mid T\in \mathcal{F}_{X/C}\}.\end{equation}
Now, we prove $F_{X/C} \cong (\mathcal{S}, \supseteq)$. Let $f: \mathcal{F}_{X/C}\rightarrow \mathcal{S}$ be defined by $f(T)=S_T$ for all $T\in \mathcal{F}_{X/C}$.
Clearly, $f$ is a surjective map by (\ref{equa 4}) and (\ref{equa 5}). Suppose that $T_1 \neq T_2$ in $\mathcal{S}$. We claim that $f(T_1) \neq f(T_2)$.
Otherwise, $f(T_1)=f(T_2)$. Then $S_{T_1}= S_{T_2}$. Without loss of generality, we suppose that $[x]\in T_1$ but
$[x]\notin T_2$ since $T_1 \neq T_2$. Thus $x\in [x] \subseteq S_{T_1}$ by (\ref{equa 4}), which means that $x\in S_{T_2}$. Hence, there exists $[y]\in T_2$
such that $x\in [y]$. Thus $[x]=[y]\in T_2$, a contradiction. So that $f$ is injective.

Therefore, $f$ is bijective.

Next, we prove that both $f$ and $f^{-1}$ preserve the order $\supseteq$.

If $T_1 \supseteq T_2$ then obviously $f(T_1)=S_{T_1}\supseteq S_{T_2}=f(T_2)$, i.e., $f$ preserves the order $\supseteq$.

Furthermore, we claim that $T_1 \supseteq T_2$ when $f(T_1)\supseteq f(T_2)$. Suppose $T_1 \nsupseteq T_2$ but $f(T_1)\supseteq f(T_2)$.
Then there exists an element $[z]\in T_2$ but $[z]\notin T_1$. Note that $z\in S_{T_2}=f(T_2)$ by formula (\ref{equa 4}). This implies that $z\in f(T_1)$ since $f(T_1)\supseteq f(T_2)$. Thus, from (\ref{equa 4}), there exists $[u]\in T_1$ such
that $z\in [u]$, and so that $[z]=[u]\in T_1$, a contradiction. Consequently $f^{-1}$ preserves the order $\supseteq$.

Thus $f$ is an isomorphism from $F_{X/C}$ to $(\mathcal{S}, \supseteq)$, i.e., \begin{equation}\label{equa 7}F_{X/C} \cong (\mathcal{S}, \supseteq).\end{equation}

In what follows, we shall prove that $F_{X/C}$ can be embedded into $F_X$, such that all infima, suprema, the top and bottom elements are preserved under the embedding.

First, by (\ref{equa 4}), we have \begin{equation}\label{eq112}\emptyset=S_{\emptyset}\in \mathcal{S}\mbox{ and }X=S_{X/C}\in \mathcal{S}.\end{equation} Now,
let $K\in \mathcal{S}$. Then, using (\ref{equa 5}), there is a $T \in \mathcal{F}_{X/C}$ such that $K=S_T$. By Lemma \ref{le 3}, we
know that for all $x, y\in X$, \begin{equation}\label{2222}y\geq x \mbox{ in }X\mbox{ implies } [y]\geq [x] \mbox{ in } X/C.\end{equation}
Suppose that $x\in K$, $y\in X$ and $y\geq x$. Then by formula (\ref{2222}), $[y]\geq [x]$, and from formula (\ref{equa 4}), the condition $x\in K=S_T$ yields that $[x]\in T$. Hence $[y]\in T$ by the definition of an up-set since $T\in \mathcal{F}_{X/C}$. Therefore, $y\in S_T=K$. This means that $K$ is an up-set of $X$, i.e., $K\in \mathcal{F}_X$.
By the arbitrariness of $K$, we
have \begin{equation}\label{equa 9}\mathcal{S}\subseteq \mathcal{F}_X.\end{equation}

Secondly, let $\{S_i\mid i\in I\}\subseteq\mathcal{S}$. From (\ref{equa 5}), there
exists $T_i \in \mathcal{F}_{X/C}$ such that $S_i=S_{T_{i}}$ for all $i\in I$. Thus $\bigcap_{i\in I}S_i=\bigcap_{i\in I}S_{T_{i}}$
and $\bigcup_{i\in I}S_i=\bigcup_{i\in I}S_{T_{i}}$. Let $W=\bigcap_{i\in I}T_i$ and $R=\bigcup_{i\in I}T_i$. Then $W, R\in \mathcal{F}_{X/C}$
since $\mathcal{F}_{X/C}$ is the set of all up-sets of $X/C$. Thus, by (\ref{equa 5}), $S_W, S_R \in \mathcal{S}$. Note that $[x]\cap [y]=\emptyset$ for any $[x]\neq [y]$. Thus, by (\ref{equa 4}), we further know that $S_W=\bigcap_{i\in I}S_{T_i}$ and $S_R=\bigcup_{i\in I}S_{T_i}$.
Therefore, \begin{equation}\label{equa 8}\bigcap_{i\in I}S_i \in \mathcal{S} \mbox{ and } \bigcup_{i\in I}S_i \in \mathcal{S}.\end{equation}

Finally, from formulas (\ref{equa 7}), (\ref{eq112}), (\ref{equa 9}) and (\ref{equa 8}), we know that $F_{X/C}$ can be embedded into $F_X$, such that all infima, suprema, and the top and bottom elements are preserved under the embedding.
\endproof

The following example will illustrate Theorem \ref{the0 3}.
\begin{example}\label{exa 1}
\emph{Let us consider the poset $(X, \leq)$ where $X=\{a, b, c, d, e\}$, and the poset $(X/C, \leq)$ represented in Fig. 2, where $C$: $X\rightarrow X$ is
a closure operator on $X$ satisfying that $C(x)=\left\{
\begin{array}{rcl}
x       &      & {x\neq d,}\\
e     &      & {x=d }
\end{array}\right.$ for all $x\in X$.}

\par\noindent\vskip50pt
 \begin{minipage}{11pc}
\setlength{\unitlength}{0.75pt}\begin{picture}(600,100)
\put(320,40){\circle{4}}\put(300,30){\makebox(0,0)[l]{\footnotesize $[a]$}}
\put(400,40){\circle{4}}\put(405,30){\makebox(0,0)[l]{\footnotesize $[b]$}}
\put(320,120){\circle{4}}\put(300,120){\makebox(0,0)[l]{\footnotesize $[c]$}}
\put(400,120){\circle{4}}\put(405,120){\makebox(0,0)[l]{\footnotesize $[e]=[d]$}}
\put(320,42){\line(0,1){76}}
\put(400,42){\line(0,1){76}}
\put(321.5,41.5){\line(1,1){77}}
\put(398.5,41.5){\line(-1,1){77}}
\put(420,80){$X/C$}
\put(160,40){\circle{4}}\put(155,30){\makebox(0,0)[l]{\footnotesize $a$}}
\put(240,80){\circle{4}}\put(245,80){\makebox(0,0)[l]{\footnotesize $d$}}
\put(240,40){\circle{4}}\put(245,30){\makebox(0,0)[l]{\footnotesize $b$}}
\put(160,120){\circle{4}}\put(150,120){\makebox(0,0)[l]{\footnotesize $c$}}
\put(240,120){\circle{4}}\put(245,120){\makebox(0,0)[l]{\footnotesize $e$}}
\put(160,42){\line(0,1){76}}
\put(240,42){\line(0,1){36}}
\put(240,82){\line(0,1){36}}
\put(161.5,41.5){\line(1,1){77}}
\put(238.5,41.5){\line(-1,1){77}}
\put(130,80){$X$}
\put(155,0){\emph{Fig.2 Hasse diagrams of $X$ and $X/C$}}
  \end{picture}
  \end{minipage}
  $$ $$
\emph{It is clear that $$\mathcal{F}_X=\{\emptyset, \{c\}, \{e\}, \{c, e\}, \{d, e\}, \{c, d, e\}, \{a, c, e\}, \{a, c, d, e\}, \{b, c, d, e\}, X\} \mbox{ and } $$
$$\mathcal{F}_{X/C}=\{\emptyset, \{[c]\}, \{[e]\}, \{[c], [e]\}, \{[a], [c], [e]\}, \{[b], [c], [e]\}, X/C\}.$$
Therefore, $F_{X/C}$ can be embedded into $F_X$, such that the top and bottom elements are preserved under the embedding.}
\end{example}

As it is well known, any finite distributive lattice $L$ can be isomorphically represented by the collection of all up-sets of $M(L)$, ordered dually to inclusion (Bikhoff's Representation Theorem). Then from Theorem \ref{the0 3} we have the following corollary.
\begin{corollary}\label{coro 2}
Let $L_1$ and $L_2$ be two finite distributive lattices. If there exists a closure operator
$C$ on $M(L_1)$ such that $M(L_1)/C \cong M(L_2)$ then $L_2$ can be embedded into $L_1$, such that the top and bottom
elements are preserved under the embedding.
\end{corollary}

The following example will show that the inverse of Corollary \ref{coro 2} does not hold generally.
\begin{example}\label{exa 2}
\emph{Let us consider again the poset $(X, \leq)$ in Example \ref{exa 1}. Let $L_1=F_X$ and $L_2=(\mathcal{U}, \supseteq)$, where
$$\mathcal{U}=\{\emptyset, \{c\}, \{c, e\}, \{c, d, e\}, \{a, c, d, e\},  X\}.$$ Clearly, $M(L_1)\cong (X, \leq)$ and the lattice $L_2$ can
be embedded into $L_1$, such that the top and bottom
elements are preserved under the embedding. However, $M(L_2)\ncong M(L_1)/C$ for any closure operator $C$ on $M(L_1)$ since $L_2$ is a six-element chain.}
\end{example}

\begin{theorem}\label{theo 4}
Let $(X, \leq)$ be a poset, let $C$ be a closure operator on $X$, let $(L, \leq)$ be a complete lattice
and let $\mathcal{S}$ be defined as (\ref{equa 5}). If there exists a closure operator $C_1$ such that $L/C_1 \cong F_X$, then $(\mathcal{F}^{u}_{L}(X/C)/\sim, \leq)\cong [0,[\mu]_\sim]$, in which, $0$ is the bottom element of the lattice $(\mathcal{F}^{u}_{L}(X)/\sim, \leq)$
and $\mu \in \mathcal{F}^{u}_{L}(X)$ satisfying $\mu_L=\mathcal{S}$.
\end{theorem}
\proof
First, by Theorem \ref{the0 2} and Proposition \ref{pro 3}, we know that $(\mathcal{F}^{u}_{L}(X)/\sim, \leq)$ is a complete lattice.

Next, we shall prove that $(\mathcal{F}^{u}_{L}(X/C)/\sim, \leq)\cong [0,[\mu]_\sim]$.

From Theorem \ref{the 1}, the condition $L/C_1 \cong F_X$ yields that there exists $h\in \mathcal{F}^{u}_{L}(X)$ such that $h_L=\mathcal{F}_X$. Then, using formula (\ref{equa 9}), $\mathcal{S}\subseteq h_L$. Moreover, by (\ref{eq112}) and (\ref{equa 8}), it follows from Theorem \ref{the0 1} that there exists $\mu \in \mathcal{F}^{u}_{L}(X)$ such
that \begin{equation}\label{ee1}\mu_L=\mathcal{S}.\end{equation} Hence, by Theorem \ref{the 1}, there exists a closure operator $C_0$ such that $L/C_0 \cong (\mathcal{S},\supseteq)$. Furthermore, $L/C_0 \cong F_{X/C}$ by (\ref{equa 7}). Then by Theorem \ref{the0 2} and Proposition \ref{pro 3}, we know that $(\mathcal{F}^{u}_{L}(X/C)/\sim, \leq)$ is a complete lattice.
Again by Theorem \ref{the 1}, there exists $g\in \mathcal{F}^{u}_{L}(X/C)$ such that \begin{equation}\label{ee101}g_L=\mathcal{F}_{X/C}.\end{equation}

Denote \begin{equation*}\mathcal{S}^{*}=\{T\subseteq \mathcal{S}\mid T \mbox{ is closed under intersections and } X\in T\}.\end{equation*}
Thus, by Theorem \ref{the0 1} and (\ref{ee1}), if $T\in \mathcal{S}^{*}$ then there exists $\nu\in \mathcal{F}^{u}_{L}(X)$ such that \begin{equation}\label{ee111}\nu_L=T.\end{equation}
Formulas (\ref{ee111}) and (\ref{ee1}) imply that $\nu_L \subseteq \mu_L$, so that $[\nu]_\sim\leq [\mu]_\sim$, i.e., $[\nu]_\sim\in [0, [\mu]_\sim]$. Note that, if $[f]_\sim\in [0, [\mu]_\sim]$ then $f_L \in \mathcal{S}^{*}$. Hence
\begin{equation}\label{ee2}[0, [\mu]_\sim]\cong(\mathcal{S}^{*},\subseteq).\end{equation}
Using $\mathcal{F}_{X/C}$ and $[g]_\sim$ instead of $\mathcal{S}$ and $[\mu]_\sim$, respectively, then similar to the proof of (\ref{ee2}), we have that $$[0_1, [g]_\sim]\cong((\mathcal{F}_{X/C})^{*},\subseteq)$$
in which $0_1$ is the bottom element of $(\mathcal{F}^{u}_{L}(X/C)/\sim, \leq)$. Because $[g]_\sim$ is the top element of $(\mathcal{F}^{u}_{L}(X/C)/\sim, \leq)$ by (\ref{ee101}) and Theorem \ref{the 2}, we have $(\mathcal{F}^{u}_{L}(X/C)/\sim, \leq)=[0_1,[g]_\sim]$. Therefore, \begin{equation}\label{ee3}(\mathcal{F}^{u}_{L}(X/C)/\sim, \leq)\cong((\mathcal{F}_{X/C})^{*},\subseteq).\end{equation}
Formula (\ref{equa 7}) implies that $((\mathcal{F}_{X/C})^{*},\subseteq)\cong (\mathcal{S}^{*},\subseteq)$. Consequently, we have that $$(\mathcal{F}^{u}_{L}(X/C)/\sim, \leq)\cong [0,[\mu]_\sim]$$
by formulas (\ref{ee2}) and (\ref{ee3}).
\endproof

\section{Conclusions}
This paper gave a necessary and sufficient condition under which $(\mathcal{F}^{u}_{L}(X)/\sim, \leq)$ is a complete lattice by terminologies of closure operators, and proved that $(\mathcal{F}^{u}_{L}(X/C)/\sim, \leq)$ is isomorphic to an interval generated by an $L$-fuzzy up-set where $C$ is a closure operator on the poset $X$. The last result shows us that for a given closure operator $C$ on a poset $X$ if there exists a closure operator $C_1$ on a complete lattice $L$ such that $L/C_1 \cong F_X$, then there is a $\mu \in \mathcal{F}^{u}_{L}(X)$ such that $(\mathcal{F}^{u}_{L}(X/C)/\sim, \leq)\cong [0,[\mu]_\sim]$. However, there exists also an $L$-fuzzy up-set $\mu\in \mathcal{F}^{u}_{L}(X)$ such that $$(\mathcal{F}^{u}_{L}(X/C)/\sim, \leq)\ncong [0,[\mu]_\sim]$$
for all closure operators $C$ on the poset $X$.

For instance, consider again Example \ref{exa 2}. Let $L$ be a complete lattice, and let $C_1$ be a closure operator on $L$ such that $L/C_1\cong F_X$
and $\mu \in \mathcal{F}^{u}_{L}(X)$ with $\mu_L=\mathcal{U}$. Then by Example \ref{exa 2}, $L_2$ is not isomorphic to $F_{X/C}$, so that for any closure operator $C$ on the poset
$X$, $(\mathcal{F}^{u}_{L}(X/C)/\sim, \leq)\ncong [0,[\mu]_\sim]$.

Therefore, an interesting problem is whether we can weaken the conditions of closure operators such that for a given $L$-fuzzy up-set $\mu\in \mathcal{F}^{u}_{L}(X)$ there exists a weakening closure operator, say $\Delta$, such that $$(\mathcal{F}^{u}_{L}(X/\Delta)/\sim, \leq)\cong [0,[\mu]_\sim].$$


\begin{thebibliography}{10}
\bibitem{Borchar} B. Borchardt, A. Maletti, B. \v{S}e\v{s}elja, A. Tepav\v{c}evi\'{c}, H. Vogler, Cut sets as recongnizable
tree languages, Fuzzy Sets and Systems 157 (2006) 1560-1571.
\bibitem{Chechik} M. Chechik, B. Devereux, S. Easterbrook, A. Y. C. Lai, V. Petrovykh, Efficient multiple-valued model-checking
using lattice representations, in: Lecture Notes in Computer Scince, Vol. 2154, Springer, Berlin, 2001, pp. 451-465.

\bibitem{Gratze} B. A. Davey, H. A. Priestley, Introrduction to Lattices and Order, Cambridge University Press, Cambridge, 2012.
\bibitem{Goguen} J. A. Goguen, L-fuzzy sets, J. Math. Anal. Appl. 18 (1967) 145-174.



\bibitem{Gorjan} M. Gorjanac-Ranitovi\'{c}, A. Tepav\v{c}evi\'{c}, General form of lattice-valued fuzzy sets under the cutworthy
approach, Fuzzy Sets and Systems 158 (2007) 1213-1216.
\bibitem{He} Peng He, Xue-ping Wang, On the uniqueness of $L$-fuzzy sets in the representation of families of sets, submitted (arXiv:1610.03532).

\bibitem{Jorge} J. Jim\'{e}nez, S. Montes, B. \v{S}e\v{s}elja, A. Tepav\v{c}evi\'{c}, On lattice valued up-sets
and down sets, Fuzzy Sets and Systems. 161 (2010) 1699-1710.

\bibitem{Seselja5} B. \v{S}e\v{s}elja, A. Tepav\v{c}evi\'{c}, A note on natural equivalence relation on fuzzy power set, Fuzzy Sets and Systems. 148 (2004) 201-210.

\end{thebibliography}
\end{document}